\documentclass{article}
\usepackage{graphicx}
\usepackage{amsmath}
\usepackage{amsfonts}
\usepackage{amssymb}
\newtheorem{theorem}{Theorem}

\newtheorem{corollary}[theorem]{Corollary}

\newtheorem{lemma}[theorem]{Lemma}

\newtheorem{proposition}[theorem]{Proposition}

\begin{document}

\title{A general Pl\"{u}cker formula}
\author{Naichung Conan Leung}
\maketitle
\begin{abstract}
We prove a formula which compares intersection numbers of conormal varieties
of two projective varieties and their dual varieties.

When one of them is linear, we can recover the usual Pl\"{u}cker formula for
the degree of the dual variety.

The basic strategy of the proof is to study a category of Lagrangian
subvarieties in the cotangent bundle of a projective space under a birational transformation.
\end{abstract}

\newpage

\section{The formula}

The conormal variety $C_{S}$ of any subvariety $S$ in $\mathbb{P}^{n}$ is a
Lagrangian subvariety in\ $T^{\ast}\mathbb{P}^{n}$ with respect to the
canonical holomorphic symplectic form. If $S$ is smooth then its Euler
characteristic $\chi\left(  S\right)  $ equals to the intersection number
$C_{S}\cdot\mathbb{P}^{n}$ up to a sign $\left(  -1\right)  ^{\dim S}$. In
general $\chi\left(  S\right)  $ is replaced by the Euler characteristic of
the Euler obstruction $Eu_{S}$ defined by MacPherson \cite{Ma}, we denote it
as $\bar{\chi}\left(  S\right)  $. More generally $C_{S_{1}}\cdot C_{S_{2}}$
is well-defined and equals to the Euler characteristic of the intersection up
to a sign provided that $S_{1}$ and $S_{2}$ intersect transversely along a
smooth subvariety in $\mathbb{P}^{n}$.

In this paper we prove the following formula.

\medskip

\begin{theorem}
\label{1ThmMain}Suppose $S_{1}$ and $S_{2}$ are two subvarieties in
$\mathbb{P}^{n}$ which intersect transversely and the same holds true for
their dual varieties $S_{1}^{\vee}$ and $S_{2}^{\vee}$ in $\mathbb{P}^{n\ast}%
$. Then we have\medskip\bigskip\
\[
C_{S_{1}}\cdot C_{S_{2}}+\frac{\left(  C_{S_{1}}\cdot\mathbb{P}^{n}\right)
\left(  C_{S_{2}}\cdot\mathbb{P}^{n}\right)  }{\left(  -1\right)
^{n+1}\left(  n+1\right)  }=C_{S_{1}^{\vee}}\cdot C_{S_{2}^{\vee}}%
+\frac{\left(  C_{S_{1}^{\vee}}\cdot\mathbb{P}^{n\ast}\right)  \left(
C_{S_{2}^{\vee}}\cdot\mathbb{P}^{n\ast}\right)  }{\left(  -1\right)
^{n+1}\left(  n+1\right)  }.
\]
\end{theorem}

\medskip

\bigskip

This formula arises as we studied in \cite{Le} the Legendre transformation on
the category of Lagrangian subvarieties in a hyperk\"{a}hler manifold under a
birational transformation of the ambient manifold.

\bigskip

By applying this formula to cases when $S_{2}$ are linear subspaces in
$\mathbb{P}^{n}$ with different dimensions, we obtain the following three
corollaries which determine the \textit{degree}, the \textit{Euler
characteristic} and the \textit{dimension} of the dual variety. We note that
the degree of $S_{1}$ equals $C_{S_{1}}\cdot C_{S_{2}}$ with $S_{2}$ a linear
subspace of complementary dimension.

First we recover the generalized Pl\"{u}cker formula of Parusi\'{n}ski
\cite{Pa} and Ernstr\"{o}m \cite{Er}.

\begin{corollary}
For any $k$ we have
\[
\deg S^{\vee}=\left(  -1\right)  ^{n-k+1}\left(  k\bar{\chi}\left(  S\right)
-\left(  k+1\right)  \bar{\chi}\left(  S^{1}\right)  +\bar{\chi}\left(
S^{k+1}\right)  \right)  ,
\]
where $S^{k}$ is the intersection of $S$ with a generic codimensional $k$
linear subspace.
\end{corollary}

When $S$ is a plane curve in $\mathbb{P}^{2}$ this gives the classical
Pl\"{u}cker formula for the degree of the dual curve $S^{\vee}$,
\[
\deg S^{\vee}=d\left(  d-1\right)  -2\delta-3\kappa,
\]
where $d,\delta,\kappa$ denote the degree, the number of nodes, the number of
cusps of $S$. This classical formula has another generalization for higher
dimensional $S$ with only isolated singularities by Teissier \cite{Te} in the
hypersurface case and by Kleiman \cite{Kl} in general. For subjects closely
related to the Pl\"{u}cker formula, readers can consult \cite{GKZ} and
\cite{Kl2}.

\bigskip

Second, we can derive the Euler characteristic of the dual variety.

\begin{corollary}
For any subvariety $S$ in $\mathbb{P}^{n}$ we have
\[
\bar{\chi}\left(  S^{\vee}\right)  =n\bar{\chi}\left(  S\right)  -\left(
n+1\right)  \bar{\chi}\left(  S^{1}\right)  \text{.}%
\]
\end{corollary}

Third we can determine the dimension of the dual variety by comparing
$\chi\left(  S^{k}\right)  $ with a linear function in $k$. For example, in
the hypersurface case we have,

\begin{corollary}
Let $S$ be a hypersurface in $\mathbb{P}^{n}$, then $S^{\vee}$ has codimension
$c$ if and only if for any $k\leq c$ we have
\[
\bar{\chi}\left(  S^{k}\right)  =k\bar{\chi}\left(  S^{1}\right)  +\left(
1-k\right)  \bar{\chi}\left(  S\right)
\]
and it becomes a strict inequality when $k=c+1$.
\end{corollary}

There are dimension formulas for dual varieties in terms of Hessian matrices
given by Segre \cite{Se} and Katz \cite{Ka}.

\bigskip

When $S_{2}$ is a smooth quadric hypersurface, its dual variety $S_{2}^{\vee}
$ is again a smooth quadric hypersurface. In this case we have the following corollary.

\begin{corollary}
Suppose a subvariety $S$ in $\mathbb{P}^{n}$ and its dual variety $S^{\vee}$
are both hypersurfaces then
\[
\bar{\chi}\left(  S\right)  -\hat{\chi}\left(  S\cap Q\right)  \left(
1+\frac{1+\left(  -1\right)  ^{n}}{2n}\right)  =\bar{\chi}\left(  S^{\vee
}\right)  -\hat{\chi}\left(  S^{\vee}\cap Q^{\prime}\right)  \left(
1+\frac{1+\left(  -1\right)  ^{n}}{2n}\right)  ,
\]
where $Q,Q^{\prime}$ are general quadric hypersurfaces and $\hat{\chi}\left(
S\cap Q\right)  =\left(  -1\right)  ^{n-2}C_{S}\cdot C_{Q}$.
\end{corollary}

Remark: We expect that $\hat{\chi}\left(  S\cap Q\right)  $ is simply equal to
$\bar{\chi}\left(  S\cap Q\right)  $.

Proof: By applying our formula to the case when $S_{1}$ is a hyperplane and
$S_{2}=Q_{n-1}$ is a general quadric hypersurface in $\mathbb{P}^{n}$, we
obtain
\begin{align*}
\chi\left(  Q_{n}\right)   &  =n+1\text{ if }n\text{ is odd,}\\
\chi\left(  Q_{n}\right)   &  =n+2\text{ if }n\text{ is even.}%
\end{align*}
This can also be obtained by using explicit descriptions of $Q_{n}$'s. We
apply our formula again by replacing $S_{1}$ with the hypersurface $S$ and we
obtain the result. $\blacksquare$

\newpage

\section{Proof of the formula}

The basic idea of our proof is to study the intersection of two Lagrangian
subvarieties $C_{1}$ and $C_{2}$ in the holomorphic symplectic manifold
$T^{\ast}\mathbb{P}^{n}$ and their behaviors under the canonical birational
transformation from $T^{\ast}\mathbb{P}^{n}$ to $T^{\ast}\mathbb{P}^{n\ast}$:

Recall that we can write $\mathbb{P}^{n}=\left(  V\backslash0\right)
/\mathbb{C}^{\times}$ and $\mathbb{P}^{n\ast}=\left(  V^{\ast}\backslash
0\right)  /\mathbb{C}^{\times}$ with $V$ a vector space of dimension $n+1$.
There is a similar description for their cotangent bundles, namely
\begin{align*}
T^{\ast}\mathbb{P}^{n}  &  =\left\{  \left(  x,\xi\right)  \in\left(
V\backslash0\right)  \times V^{\ast}:\xi\left(  x\right)  =0\right\}
/\mathbb{C}^{\times},\\
T^{\ast}\mathbb{P}^{n\ast}  &  =\left\{  \left(  x,\xi\right)  \in
V\times\left(  V^{\ast}\backslash0\right)  :\xi\left(  x\right)  =0\right\}
/\mathbb{C}^{\times}.
\end{align*}
The zero section $P$ in $T^{\ast}\mathbb{P}^{n}$ (resp. $P^{\ast}$ in
$T^{\ast}\mathbb{P}^{n\ast}$) is given by $\left\{  \xi=0\right\}  $ (resp.
$\left\{  x=0\right\}  $). The identity homomorphism on $V\times V^{\ast}$
descends to a birational map,
\[
\Phi:T^{\ast}\mathbb{P}^{n}\dashrightarrow T^{\ast}\mathbb{P}^{n\ast}\text{,}%
\]
which is biregular from $T^{\ast}\mathbb{P}^{n}\backslash P$ to $T^{\ast
}\mathbb{P}^{n\ast}\backslash P^{\ast}$.

\bigskip

Even though $T^{\ast}\mathbb{P}^{n}$ is incomplete, the intersection of
$C_{S_{1}}$ and $C_{S_{2}}$ only occurs along the zero section $P$ in
$T^{\ast}\mathbb{P}^{n}$ when $S_{1}$ and $S_{2}$ intersect transversely in
$\mathbb{P}^{n}$.\footnote{For now on we will simply write $C_{i}$ for
$C_{S_{i}}$, the conormal variety of $S_{i}$.} This holds true even after we
compactify $T^{\ast}\mathbb{P}^{n}$ to $M$ as follow,
\[
M=\mathbb{P}\left(  T^{\ast}\mathbb{P}^{n}\oplus O_{\mathbb{P}^{n}}\right)
\text{.}%
\]
If we denote the closure of $C_{i}$'s in $M$ as $\bar{C}_{i}$'s, then we have
\[%
\begin{array}
[c]{ccccc}%
\bar{C}_{1}\cap\bar{C}_{2} & = & C_{1}\cap C_{2} & = & S_{1}\cap S_{2}\\
\cap &  & \cap &  & \cap\\
M & \supset &  T^{\ast}\mathbb{P}^{n} & \supset & \mathbb{P}^{n}.
\end{array}
\]
Therefore we can simply write $\bar{C}_{1}\cdot\bar{C}_{2}$ as $C_{1}\cdot
C_{2}$. Similarly we have $C_{i}\cdot P=\bar{C}_{i}\cdot P$.

\bigskip

Remark: We recall the following useful result (see e.g. \cite{Le}): If $C_{1}
$ and $C_{2}$ are Lagrangian subvarieties in $T^{\ast}\mathbb{P}^{n}$ which
intersect cleanly. Suppose their closures $\bar{C}_{1}$ and $\bar{C}_{2}$ have
the same intersection, i.e. $C_{1}\cap C_{2}=\bar{C}_{1}\cap\bar{C}_{2} $,
then
\[
C_{1}\cdot C_{2}=\left(  -1\right)  ^{\dim C_{1}\cap C_{2}}\chi\left(
C_{1}\cap C_{2}\right)  \text{.}%
\]
In particular we have $P\cdot P=\left(  -1\right)  ^{n}\left(  n+1\right)  $.

\bigskip

It is not difficult to see that the above birational map $\Phi$ on $T^{\ast
}\mathbb{P}^{n}$ extends to a birational map between two complete varieties $M
$ and $M^{\prime}$:
\[
\Phi_{M}:M\dashrightarrow M^{\prime}\text{,}%
\]
which is biregular outside $P$ and $P^{\ast}$.

Note that $M^{\prime}$ can also be described as a flop as follow:
\[%
\begin{array}
[c]{ccccc}%
M & \overset{\pi}{\leftarrow} & \widetilde{M} & \overset{\pi^{\prime}%
}{\rightarrow} & M^{\prime}\\
\cup &  & \,\cup j &  & \cup\\
P & \overset{p}{\leftarrow} & \widetilde{P} & \overset{p^{\prime}}%
{\rightarrow} & P^{\ast}.
\end{array}
\]
Here $\widetilde{M}$ is the blow up of $M$ along $P$. The exceptional locus
$\widetilde{P}=\pi^{-1}\left(  P\right)  $ is a $\mathbb{P}^{n-1}$-bundle over
$P$. It admits another $\mathbb{P}^{n-1}$-bundle structure over the dual
projective space $P^{\ast}$ and $\pi^{\prime}$ is the blow down of
$\widetilde{M}$ along this second fiber structure on $\widetilde{P}$.

In the next section, we will show that the Legendre functor $\mathbf{L}%
=\pi_{\ast}^{\prime}\circ\pi^{-1}$ on the derived categories of coherent
sheaves on $M$ and $M^{\prime}$ is an equivalence of categories. This implies
that
\[
Ext_{O_{M}}^{k}\left(  S_{1},S_{2}\right)  \cong Ext_{O_{M^{\prime}}}%
^{k}\left(  \mathbf{L}\left(  S_{1}\right)  ,\mathbf{L}\left(  S_{2}\right)
\right)  ,
\]
for any coherent sheaves $S_{1},S_{2}$ on $M$. In order to use this
equivalence to prove the theorem \ref{1ThmMain}, we need to compute the Chern
characters of $\mathbf{L}\left(  O_{P}\right)  $ and $\mathbf{L}\left(
O_{C_{i}}\right)  $'s. First we relate the Lagrangian intersection number
$C_{1}\cdot C_{2}$ with the category of coherent sheaves on $M$.

\begin{lemma}
For any $n$ dimensional subvarieties $\bar{C}_{1},\bar{C}_{2}$ in $M$ we have
\[%
{\textstyle\sum_{k}}
\dim\left(  -1\right)  ^{k}Ext_{O_{M}}^{k}\left(  O_{\bar{C}_{1}},O_{\bar
{C}_{2}}\right)  =\left(  -1\right)  ^{n}\bar{C}_{1}\cdot\bar{C}_{2}.
\]
\end{lemma}

Proof: We recall the Riemann-Roch formula for the global Ext groups: For any
coherent sheaves $S_{1}$ and $S_{2}$ on $M$ we have,
\[%
{\textstyle\sum_{k}}
\dim\left(  -1\right)  ^{k}Ext_{O_{M}}^{k}\left(  S_{1},S_{2}\right)
=\int_{M}\overline{ch}\left(  S_{1}\right)  ch\left(  S_{2}\right)  Td_{M}%
\]
where $\overline{ch}\left(  S_{1}\right)  =\Sigma\left(  -1\right)  ^{k}%
ch_{k}\left(  S_{1}\right)  $. For $S_{i}=O_{\bar{C}_{i}}$ the structure sheaf
of a subvariety $\bar{C}_{i}$ of dimension $n$, we have
\begin{align*}
ch_{k}\left(  O_{\bar{C}_{i}}\right)   &  =0\text{ for }k<n,\\
ch_{n}\left(  O_{\bar{C}_{i}}\right)   &  =\left[  \bar{C}_{i}\right]
\text{,}%
\end{align*}
where $\left[  \bar{C}_{i}\right]  $ denotes the Poincar\'{e} dual of the
variety $\bar{C}_{i}$.

Therefore
\begin{align*}
&  \dim\left(  -1\right)  ^{k}Ext_{O_{M}}^{k}\left(  O_{\bar{C}_{1}}%
,O_{\bar{C}_{2}}\right) \\
&  =\int_{M}\overline{ch}\left(  O_{\bar{C}_{1}}\right)  ch\left(  O_{\bar
{C}_{2}}\right)  Td_{M}\\
&  =\int_{M}\left(  \left(  -1\right)  ^{n}\left[  \bar{C}_{1}\right]
+h.o.t.\right)  \left(  \left[  \bar{C}_{2}\right]  +h.o.t.\right)  \left(
1+h.o.t.\right) \\
&  =\left(  -1\right)  ^{n}\int_{M}\left[  \bar{C}_{1}\right]  \cup\left[
\bar{C}_{2}\right]  =\left(  -1\right)  ^{n}\bar{C}_{1}\cdot\bar{C}%
_{2}\text{.}%
\end{align*}
Here $h.o.t.$ refers to \textit{higher order terms} which do not contribute to
the outcome of the integral. Hence the result. $\blacksquare$

\begin{lemma}%
\[
ch\left(  \mathbf{L}\left(  O_{P}\right)  \right)  =\pm\left[  P^{\ast
}\right]  +h.o.t.\text{.}%
\]
\end{lemma}

Proof: From the previous lemma, we have
\begin{align*}
&
{\textstyle\sum}
\left(  -1\right)  ^{k}\dim Ext_{O_{M}}^{k}\left(  O_{P},O_{P}\right) \\
&  =\left(  -1\right)  ^{n}P\cdot P\\
&  =\chi\left(  P\right) \\
&  =\left(  n+1\right)  \text{.}%
\end{align*}

On the other hand, the support of $\mathbf{L}\left(  O_{P}\right)  $ must be
inside $P^{\ast}\subset M^{\prime}$ because $M$ and $M^{\prime}$ are
isomorphic outside their exceptional loci $P$ and $P^{\ast}$. This implies
that $ch\left(  \mathbf{L}\left(  O_{P}\right)  \right)  =\alpha\left[
P^{\ast}\right]  +h.o.t.$ for some integer $\alpha$. Therefore,%

\begin{align*}
&
{\textstyle\sum_{k}}
\left(  -1\right)  ^{k}\dim Ext_{O_{M^{\prime}}}^{k}\left(  \mathbf{L}\left(
O_{P}\right)  ,\mathbf{L}\left(  O_{P}\right)  \right) \\
&  =\int_{M^{\prime}}\overline{ch}\left(  \mathbf{L}\left(  O_{P}\right)
\right)  ch\left(  \mathbf{L}\left(  O_{P}\right)  \right)  Td_{M}\\
&  =\int_{M^{\prime}}\left(  -1\right)  ^{n}\alpha\left[  P^{\ast}\right]
\cdot\alpha\left[  P^{\ast}\right] \\
&  =\left(  n+1\right)  \alpha^{2}\text{.}%
\end{align*}
Now $\mathbf{L}$ being an equivalence of categories implies that
\[%
{\textstyle\sum}
\left(  -1\right)  ^{k}\dim Ext_{O_{M^{\prime}}}^{k}\left(  \mathbf{L}\left(
O_{P}\right)  ,\mathbf{L}\left(  O_{P}\right)  \right)  =%
{\textstyle\sum}
\left(  -1\right)  ^{k}\dim Ext_{O_{M}}^{k}\left(  O_{P},O_{P}\right)  .
\]

This forces $\alpha$ to be $\pm1$. Hence the claim. $\blacksquare$

\begin{lemma}
Suppose $C$ is a $n$ dimensional irreducible subvariety in $M$ not containing
$P$. We consider
\[
S=%
{\textstyle\bigoplus^{n+1}}
O_{C}-\left(  -1\right)  ^{n}%
{\textstyle\bigoplus^{C\cdot P}}
O_{P}\text{.}%
\]
Then
\[
ch\left(  \mathbf{L}\left(  S\right)  \right)  =\left(  n+1\right)  \left[
C^{\vee}\right]  -\left(  -1\right)  ^{n}\left(  C^{\vee}\cdot P^{\ast
}\right)  \left[  P^{\ast}\right]  +h.o.t.
\]
\end{lemma}

Proof: First we have
\begin{align*}
&  \sum\left(  -1\right)  ^{k}\dim Ext_{O_{M}}^{k}\left(  S,O_{P}\right) \\
&  =\left(  -1\right)  ^{n}\left[  \left(  n+1\right)  C\cdot P-\left(
-1\right)  ^{n}\left(  C\cdot P\right)  \left(  P\cdot P\right)  \right] \\
&  =0\text{.}%
\end{align*}
Therefore, using the equivalence of categories induced by $\mathbf{L}$ again,
we have
\[
\sum\left(  -1\right)  ^{k}\dim Ext_{O_{M^{\prime}}}^{k}\left(  \mathbf{L}%
\left(  S\right)  ,\mathbf{L}\left(  O_{P}\right)  \right)  =0\text{.}%
\]
Away from $P^{\ast}$, the support of $\mathbf{L}\left(  S\right)  $ is clearly
$C^{\vee}$ with multiplicity $n+1$, therefore we have
\[
ch\left(  \mathbf{L}\left(  S\right)  \right)  =\left(  n+1\right)  \left[
C^{\vee}\right]  +\beta\left[  P^{\ast}\right]  +h.o.t.
\]
for some integer $\beta$. On the other hand
\begin{align*}
&  \sum\left(  -1\right)  ^{k}\dim Ext_{O_{M^{\prime}}}^{k}\left(
\mathbf{L}\left(  S\right)  ,\mathbf{L}\left(  O_{P}\right)  \right) \\
&  =\int_{M^{\prime}}\overline{ch}\left(  \mathbf{L}\left(  S\right)  \right)
ch\left(  \mathbf{L}\left(  O_{P}\right)  \right)  Td_{M^{\prime}}\\
&  =\left(  -1\right)  ^{n}\left(  \left(  n+1\right)  \left[  C^{\vee
}\right]  +\beta\left[  P^{\ast}\right]  \right)  \cdot\left(  \pm\left[
P^{\ast}\right]  \right)  \text{.}%
\end{align*}
This being zero implies that $\beta=-\left(  -1\right)  ^{n}\left(  C^{\vee
}\cdot P^{\ast}\right)  $. Therefore
\[
ch\left(  \mathbf{L}\left(  S\right)  \right)  =\left(  n+1\right)  \left[
C^{\vee}\right]  -\left(  -1\right)  ^{n}\left(  C^{\vee}\cdot P^{\ast
}\right)  \left[  P^{\ast}\right]  +h.o.t.
\]
Hence the result. $\blacksquare$

\bigskip

Now we suppose that $C_{1}$ and $C_{2}$ are two $n$ dimensional irreducible
subvarieties in $M$ which do not include $P$. We denote the corresponding
sheaves constructed as above by $S_{1}$ and $S_{2}$ respectively. We have
\begin{align*}
&  \sum\left(  -1\right)  ^{k}\dim Ext_{O_{M}}^{k}\left(  S_{1},S_{2}\right)
\\
&  =\int_{M}\overline{ch}\left(  S_{1}\right)  ch\left(  S_{2}\right)
Td_{M}\\
&  =\left(  n+1\right)  ^{2}\left(  -1\right)  ^{n}\left(  C_{1}+\frac{\left(
C_{1}\cdot P\right)  }{\left(  -1\right)  ^{n+1}\left(  n+1\right)  }P\right)
\cdot\left(  C_{2}+\frac{\left(  C_{2}\cdot P\right)  }{\left(  -1\right)
^{n+1}\left(  n+1\right)  }P\right) \\
&  =\left(  n+1\right)  ^{2}\left(  -1\right)  ^{n}\left[  C_{1}\cdot
C_{2}+\frac{\left(  -1\right)  ^{n+1}}{n+1}\left(  C_{1}\cdot P\right)
\left(  C_{2}\cdot P\right)  \right]  .
\end{align*}

By the above computations of the Chern characters of $ch\left(  \mathbf{L}%
\left(  S_{i}\right)  \right)  $'s we have a similar formula for $\sum\left(
-1\right)  ^{k}\dim Ext_{O_{M^{\prime}}}^{k}\left(  \mathbf{L}\left(
S_{1}\right)  ,\mathbf{L}\left(  S_{2}\right)  \right)  $, i.e.
\begin{align*}
&  \sum\left(  -1\right)  ^{k}\dim Ext_{O_{M^{\prime}}}^{k}\left(
\mathbf{L}\left(  S_{1}\right)  ,\mathbf{L}\left(  S_{2}\right)  \right) \\
&  =\left(  n+1\right)  ^{2}\left(  -1\right)  ^{n}\left[  C_{1}^{\vee}\cdot
C_{2}^{\vee}+\frac{\left(  -1\right)  ^{n+1}}{n+1}\left(  C_{1}^{\vee}\cdot
P^{\ast}\right)  \left(  C_{2}^{\vee}\cdot P^{\ast}\right)  \right]  .
\end{align*}
Finally using the equivalence of categories induced by $\mathbf{L}$ we have
\[
\sum\left(  -1\right)  ^{k}\dim Ext_{O_{M}}^{k}\left(  S_{1},S_{2}\right)
=\sum\left(  -1\right)  ^{k}\dim Ext_{O_{M^{\prime}}}^{k}\left(
\mathbf{L}\left(  S_{1}\right)  ,\mathbf{L}\left(  S_{2}\right)  \right)  .
\]
Combining these and we obtain
\[
C_{S_{1}}\cdot C_{S_{2}}+\frac{\left(  C_{S_{1}}\cdot P\right)  \left(
C_{S_{2}}\cdot P\right)  }{\left(  -1\right)  ^{n+1}\left(  n+1\right)
}=C_{S_{1}^{\vee}}\cdot C_{S_{2}^{\vee}}+\frac{\left(  C_{S_{1}^{\vee}}\cdot
P^{\ast}\right)  \left(  C_{S_{2}^{\vee}}\cdot P^{\ast}\right)  }{\left(
-1\right)  ^{n+1}\left(  n+1\right)  }.
\]
Thus we have completed the proof of the main theorem \ref{1ThmMain} assuming
that $\mathbf{L}$\textbf{\ }is an equivlance of categories. $\blacksquare$

\newpage

\section{Equivalence of derived categories}

In this section we adapt Bondal and Orlov arguments \cite{BO} to prove that
$D^{b}\left(  M\right)  $ and $D^{b}\left(  M^{\prime}\right)  $ are
equivalent categories. Recall that $\widetilde{M}$ is the blow up of $M$ along
$P\cong\mathbb{P}^{n}$, or the blow up of $M^{\prime}$ along $P^{\ast}%
\cong\mathbb{P}^{n\ast}$. We recall the following commutative diagram:
\[%
\begin{array}
[c]{ccccc}%
M & \overset{\pi}{\leftarrow} & \widetilde{M} & \overset{\pi^{\prime}%
}{\rightarrow} & M^{\prime}\\
\cup &  & \,\cup j &  & \cup\\
P & \overset{p}{\leftarrow} & \widetilde{P} & \overset{p^{\prime}}%
{\rightarrow} & P^{\ast}.
\end{array}
\]

\begin{proposition}
In the above situation, we have an equivalence of derived categories,
\[
\mathbf{L}=\pi_{\ast}\pi^{\prime\ast}:D^{b}\left(  M^{\prime}\right)
\rightarrow D^{b}\left(  M\right)  \text{.}%
\]
\end{proposition}

Note that $\widetilde{P}$ is a $\mathbb{P}^{n-1}$-bundle over the projective
space $P$ (or $P^{\ast}$). Therefore $Pic\left(  \widetilde{P}\right)
=\mathbb{Z}+\mathbb{Z}$. Every line bundle on $\widetilde{P}$ is isomorphic to
$p^{\ast}O_{P}\left(  a\right)  \otimes p^{\prime\ast}O_{P^{\ast}}\left(
b\right)  $ for some integers $a$ and $b$, we denote it as $O_{\widetilde{P}%
}\left(  a,b\right)  $.

We claim that $O_{\widetilde{M}}\left(  \widetilde{P}\right)  |_{\widetilde
{P}}$ is isomorphic to $O_{\widetilde{P}}\left(  -1,-1\right)  $. By symmetry
of $M$ and $M^{\prime}$ we know it must be of the form $O_{\widetilde{P}%
}\left(  a,a\right)  $ for some $a$. For the restriction of $O_{\widetilde{M}%
}\left(  \widetilde{P}\right)  $ to a fiber of $p$ equals $O_{\mathbb{P}%
^{n-1}}\left(  -1\right)  $ because $\widetilde{M}$ is obtained from $M$ by
blowing up the smooth center $P$. This implies that $a=-1$.

As a blown up manifold we have
\[
\omega_{\widetilde{M}}=\pi^{\ast}\omega_{M}\otimes O_{\widetilde{M}}\left(
\left(  n-1\right)  \widetilde{P}\right)  .
\]
This implies that
\[
\omega_{\widetilde{M}}|_{\widetilde{P}}=O_{\widetilde{P}}\left(
1-n,1-n\right)  .
\]

\bigskip

\textit{Proof of proposition:} The arguments presented here are basically the
same as in section 3 of \cite{BO}. For any $A,B\in D^{b}\left(  M^{\prime
}\right)  $ we want to show that
\begin{align*}
Hom\left(  A,B\right)   &  \cong Hom\left(  \pi_{\ast}\pi^{\prime\ast}%
A,\pi_{\ast}\pi^{\prime\ast}B\right) \\
&  \cong Hom\left(  \pi^{\ast}\pi_{\ast}\pi^{\prime\ast}A,\pi^{\prime\ast
}B\right)  \text{.}%
\end{align*}
On the other hand $Hom\left(  A,B\right)  \cong Hom\left(  \pi^{\prime\ast
}A,\pi^{\prime\ast}A\right)  $ because the functor $\pi^{\prime\ast}%
:D^{b}\left(  M^{\prime}\right)  \rightarrow D^{b}\left(  \widetilde
{M}\right)  $ is full and faithful for any blow up morphism. Therefore it is
suffices to show that
\[
Hom\left(  \bar{A},\pi^{\prime\ast}B\right)  =0,
\]
where $\bar{A}$ is defined by the following exact triangle
\[
\pi^{\ast}\pi_{\ast}\pi^{\prime\ast}A\rightarrow\pi^{\prime\ast}%
A\rightarrow\bar{A}.
\]

Using $\pi^{\ast}:D^{b}\left(  M\right)  \rightarrow D^{b}\left(
\widetilde{M}\right)  $ being full and faithful we can obtain $\bar{A}\in
D\left(  M\right)  ^{\bot}\subset D\left(  \widetilde{M}\right)  $. An earlier
result of Orlov showed that
\[
D\left(  M\right)  ^{\bot}=\left\langle D\left(  P\right)  _{-n+1}%
,\cdots,D\left(  P\right)  _{-1}\right\rangle
\]
as a semiorthogonal decomposition, where $D\left(  P\right)  _{-k}$ is the
full subcategories of $D^{b}\left(  \widetilde{M}\right)  $ given by the image
of $D^{b}\left(  P\right)  $ under $j_{\ast}\left(  O_{\widetilde{P}}\left(
-k\right)  \otimes p^{\ast}\left(  \bullet\right)  \right)  $. Since $P$
isomorphic to $\mathbb{P}^{n-1}$ we have
\[
D\left(  M\right)  ^{\bot}=\left\langle j_{\ast}O_{\widetilde{P}}\left(
a,b\right)  \right\rangle _{\substack{-n+1\leq b\leq-1 \\-1\leq a-b\leq n-1
}}
\]

From this description we know that $j_{\ast}O_{\widetilde{P}}\left(
a,b\right)  $ belongs to both $D\left(  M\right)  ^{\bot}$ and $D\left(
M^{\prime}\right)  ^{\bot}$ when $-n+1\leq a,b\leq-1$. In particular
\[
Hom\left(  \bar{A},j_{\ast}O_{\widetilde{P}}\left(  a,b\right)  \right)
=0\text{ when }-n+1\leq a,b\leq-1.
\]
Because of this $\bar{A}\in D\left(  M\right)  ^{\bot}$ has to lie inside the
subcategory of $D\left(  M\right)  ^{\bot}$ generated by those $j_{\ast
}O_{\widetilde{P}}\left(  a,b\right)  $ with $a\geq0.$ This implies that
$\bar{A}\otimes\omega_{\widetilde{M}}\in D\left(  M^{\prime}\right)  ^{\bot}$
because of $\omega_{\widetilde{M}}|_{\widetilde{P}}=O_{\widetilde{P}}\left(
1-n,1-n\right)  $. That is for any $B\in D^{b}\left(  M^{\prime}\right)  $ we
have
\[
Hom\left(  \pi^{\prime\ast}B,\bar{A}\otimes\omega_{\widetilde{M}}\right)  =0
\]
which is equivalent to
\[
Hom\left(  \bar{A},\pi^{\prime\ast}B\right)  =0,
\]
by the Serre duality. Hence our proposition. $\blacksquare$

\bigskip

\textit{Acknowledgments: The author thanks D. Abramovich, T. Bridgeland, A.
Voronov for very helpful discussions. This project is partially supported by NSF/DMS-0103355.}

\textit{\bigskip}

\bigskip

{\normalsize Address: School of Mathematics, University of Minnesota,
Minneapolis, MN 55454.}

{\normalsize Email: leung@math.umn.edu.}
\end{document}